\documentclass{article}
\usepackage{graphicx} 
\usepackage{hyperref}
\usepackage{graphicx} 
\usepackage{amssymb}
\usepackage{amsopn}
\usepackage{amsmath}
\usepackage{amsfonts}
\title{An alternative proof of Gödel's first incompleteness theorem}
\author{Zuhair Al-Johar}
\date{August 2023}

\begin{document}

\maketitle

\section{Introduction}
This proof of Gödel's first incompleteness theorem doesn't require $\omega$-consistency, nor does it refer to codes of negated sentences as in Rosser's. It begins from where Gödel's usual proof ends, and stalks it till it ends proving it. 

\section{Proof Sketch}
The proof is by negation. We assume that theory $T$ is first order, consistent, complete, and yet fulfills Gödel's criteria of representing all computable functions. So generally we take $T$ to extend either $\sf PA$ or $\sf Q$; then we show that these assumptions are contradictory. As a notation $\operatorname {Proof}_T (x,y)$ is the usual Gödel provability, also $\ulcorner . \urcorner$ refers to Gödel's coding, and $^o n$ refers to $n$-many repeated application of the successor operation over zero. The negation function $``\operatorname {neg}"$ is defined after Rosser's. Those notations are present in [4] and [5]. We use $\vdash$  to stand for provability in the usual metatheoretic sense. Generally when "we have", or provability steps are written then they refer to metatheoretic provability in $T$. 

Per those notations, we know from Gödel's work that for any standard sentence $s$ if $T\vdash s$, then there is some standard natural number $n$ such that $\operatorname {Proof}_T(n,^o \ulcorner s \urcorner)$ and vise verse. Also if we have $\operatorname {Proof}_T(n,^o\ulcorner s \urcorner)$ and $n$ is standard, then $s$ is a standard sentence. \bigskip

let $s$ be a standard sentence in the language of $T$. Let $E_s(v)$, called extension of $s$, be the sentence form:$$ \forall x \forall y \, (\operatorname {Proof}_T(x,v) \land [y= \min z: \operatorname {Proof}_T(z,^o\ulcorner s \urcorner)] \to x \geq y)$$ Now, let "$D (E_s(v))$" stand for the diagonal sentence of the sentence form $E_s(v)$ produced syntactically from $E_s(v)$ in an exact manner along lines present in [5] (see section 4) \bigskip

Accordingly; if $s$ is a standard sentence in the language of $T$, then we have: $$ \vdash_T \ D(E_s(v)) \iff \forall x \forall y \, (\operatorname {Proof}_T(x,^o \ulcorner D(E_s(v)) \urcorner) \land [y= \min z: \operatorname {Proof}_T(z,^o\ulcorner s \urcorner)] \to x \geq y)$$

Now, Let $\operatorname {Dgn}(x,y)$ be an arithmetical formula such that if there is a standard sentence $s$ and  $y= \, ^o \ulcorner s \urcorner $, then $x= \, ^o\ulcorner D(E_s(v)) \urcorner  $ [The exact definition is given at section 4]\bigskip

Define:  
\begin{align*}  remote(y) \iff \exists x \exists k \exists z: & \operatorname {Dgn}(x,y) \land \\ & k=\min n: \operatorname {Proof}_T(n,y) \land \\ & z=\min m: \operatorname {Proof}_T(m,x) \land \\& z < k \end{align*}

So, from the above definitions if $T \vdash \neg D(E_s(v))$, then $T \vdash remote(^o \ulcorner s \urcorner)$, while if $T\vdash D(E_s(v))$, then $T \vdash \neg remote(^o \ulcorner s \urcorner)$. \bigskip

Lemma 0: for every standard sentence $s$ in the language of $T$, for each concrete natural $i$, we have: $$  \vdash_T  \ s \to \forall n (\operatorname {Proof}_T(n, \operatorname{neg}(^o \ulcorner s \urcorner)) \to  n \neq \ ^oi)$$

In English: if theory $T$ proves sentence $s$ then it proves every internal proof of the negation of $s$ to be non-standard. \smallskip

Proof: if there is an internal standard proof of the negation of $s$, then per Gödel's work we'd have $T \vdash \neg s$, contradicting $T \vdash s$ \bigskip

Lemma I: for $ n=0,1,2,..$ $$ \vdash_T \forall x [remote(x) \to \forall y (\operatorname{Proof}_T(y,x) \to y \neq \ ^o n)]$$

In English: all codes of proofs of a remote sentence are non-standard. \smallskip

Proof: Suppose we have a proof of a remote sentence $s$ that is coded by a standard natural (that is $x= ^o \ulcorner s \urcorner$ and we have a standard natural $y$ such that we have $\operatorname {Proof}_T(y,x)$), then per Gödel's we'd have $s$ is a standard sentence, so the diagonal sentence $D(E_s(v))$ (see section 4) of its extension is a standard sentence and we must have $T\vdash D(E_s(v))$; otherwise we'll be having $T\vdash \neg D(E_s(v))$ and so $T$ would prove the existence of a code of a proof of $D(E_s(v))$ that is strictly smaller than the minimal code of a proof of $s$, and so would be standardly coded, and accordingly we'll have $T\vdash  D(E_s(v))$, a contradiction. So, we must have $T\vdash D(E_s(v))$ and so $\min z: \operatorname{Proof}_T (z,^o\ulcorner D(E_s(v)) \urcorner) \geq \min k: \operatorname {Proof}_T(k,^o \ulcorner s \urcorner)$, thus $\neg remote(^o \ulcorner s \urcorner)$. \bigskip

Lemma II: for any standard sentence $s$, if $T \vdash \neg s$ and $T\vdash \exists n : \operatorname{Proof}_T(n,^o\ulcorner s \urcorner)$, then $T\vdash remote(^o \ulcorner s \urcorner)$ \smallskip

In English: if $T$ refutes a standard sentence and yet proves the existence of an internal proof of that sentence, then that sentence is remote. \bigskip

Proof: Since all codes of proofs of $s$ are non-standard [Lemma 0], then any proof of the diagonal of the extension of $s$, that is $D(E_s(v))$, cannot be standardly coded, because $D(E_s(v))$ is a standard sentence and it is equivalent in $T$ to having every proof of it being not strictly less than the minimal proof of $s$, and so non-standard; so we have $T \vdash \neg D(E_s(v))$, thereby proving $remote(^o\ulcorner s \urcorner)$ \bigskip 

Define: $$x =\mathfrak K \iff x = \min k : \exists y \, (remote(y) \land \operatorname {Proof}_T(k,y))$$  

So, $\mathfrak K$ is the minimal code of a proof of a remote sentence. \bigskip

From the fixed point lemma, take:  $$\sigma \iff \forall x \, (\operatorname {Proof}_T(x,^o\ulcorner \sigma \urcorner) \to x \geq \mathfrak K)$$ Informally $\sigma$ is saying "every proof of this sentence is, bigger than or equal to, the minimal proof of remote sentences". \smallskip 

Clearly we have $T \not \vdash \sigma$. So, we must have $T \vdash \neg \sigma$, then there is a proof of $^o\ulcorner \sigma \urcorner$ in $T$ that is strictly smaller than $\mathfrak K$. But $\sigma$ is a standard sentence. So we must have $remote(^o\ulcorner \sigma \urcorner)$ [Lemma II]. Contradicting the definition of $\mathfrak K$. \smallskip

To see that; have $\sigma^*$ be the sentence $ D(E_\sigma(v))$, then: $$\sigma^* \iff \forall x \forall y (\operatorname {Proof}_T(x, ^o\ulcorner \sigma^* \urcorner) \land [y= \min z: \operatorname {Proof}_T(z,^o\ulcorner \sigma \urcorner)] \to x \geq y)$$ 
 Clearly $T \not \vdash \sigma^*$, and so we have $T\vdash \neg \sigma^*$, then we'd have $remote(^o\ulcorner \sigma \urcorner)$, but from above there is a proof of $^o\ulcorner \sigma \urcorner$ in $T$ that is strictly smaller than $\mathfrak K$, a contradiction!\bigskip

Thus $\sigma$ is undecidable in $T$. QED\smallskip

 \bigskip

The idea is that if $remote(y)$, then every code of a proof of $y$ cannot be standard, since if otherwise, then we must have the minimal proof of the diagonal of its extension being $\href{https://mathoverflow.net/a/451243/95347}{bigger}$ than the minimal proof of it. So $\mathfrak K$, the minimal a proof of a remote sentence can be, must be non-standard. The rest is clear!

The reason why $\operatorname {Dgn}$ is defined in arithmetic is because there is an exact syntactical way of producing $E_s(v)$ from $s$ and of producing $D(E_s(v))$ from $E_s(v)$, so there is an arithmetical formula such that for any standard sentence $s$ we can input its Gödel code and get the Gödel code of the diagonal of its extension. We know how this behaves exactly for all standard codes of proofs, and so by definition of $remote$ sentences none of those can have proofs coded with standard naturals. Since $T$ assumed to be complete, then we can demonstrate there is at least one remote standard sentence (e.g., the Gödel sentence), then we do have a minimal of their proofs!

\section{The proof without minimality}

The above proof depends on the minimality principle which is provable by induction. Here, this version won't assume induction. \bigskip

We need to change definitions of $E_s(v)$ and $remote(x)$ to the following: \bigskip

$E_s(v)$ is the sentence form:
$$ \forall x\, (\operatorname {Proof}_T(x,v) \to \exists y:  \operatorname {Proof}_T(y,^o\ulcorner s \urcorner) \land x \geq y)$$

Define: \smallskip 
\begin{align*}  remote(y) \iff \exists x: & \operatorname {Dgn}(x,y) \land \\ & \exists z: \operatorname {Proof}_T(z,x) \land \forall k: \operatorname {Proof}_T(k,y) \to z < k \end{align*} 

Now, the same state holds as of the old definition; that is, all above-mentioned lemmas apply and in particular no proof (in $T$) of a remote sentence can be standard! \bigskip

From the fixed point lemma, take $\sigma$ to be the sentence defined by: $$ \sigma \iff \forall x (\operatorname {Proof}_T(x,^o\ulcorner \sigma \urcorner) \to \exists y: remote(y) \land \operatorname {Proof}_T(x,y)) $$ Informally $\sigma$ says "every proof of this sentence is a proof of a remote sentence".\smallskip

 Now, if $T\vdash \sigma$, then there is a standard proof of $\sigma$, contradicting its remoteness asserted on the right hand. So, by completeness of $T$, we must have $T \vdash \neg \sigma$. So, this says that there exists a proof of $^o\ulcorner \sigma \urcorner$ in $T$ that is not a proof (in $T$) of a remote sentence, so we have $\neg remote(^o\ulcorner \sigma\urcorner)$. But, $\sigma$ is a standard sentence, and accordingly we'll have $remote (^o\ulcorner \sigma \urcorner)$ [Lemma II], a contradiction!
   
\section {Defining Dgn}

Here is the proper way of defining the relation $\operatorname {Dgn}$: \bigskip

Let $h: \mathbb N \to \mathbb N $ be the function defined by:$$h(\ulcorner s \urcorner) = \ulcorner D(E_s(v)) \urcorner$$

for each sentence $s$ in the language of theory $T$, and $h(n)=0$ otherwise. \bigskip

Define $E_s(v)$ to be the sentence form: 
$$\forall x \forall y \, (\operatorname {Proof}_T(x,v) \land [y= \min z: \operatorname {Proof}_T(z, ^o\ulcorner s \urcorner)] \to x \geq y)$$
and $D(E_s(v))$ is the diagonal formula of the sentence form $E_s(v)$, defined (along lines given \href{https://en.wikipedia.org/wiki/Diagonal_lemma#Proof}{\textbf{here}})) by:
$$\forall v (\mathcal G_f (^o\ulcorner \forall v (\mathcal G_f (k, v ) \to E_s(v)) \urcorner, v ) \to E_s(v))$$
Where $f: \mathbb N \to \mathbb N$ is the function defined by: $$ f(\ulcorner \mathcal A \urcorner)= \, \ulcorner A(^o\ulcorner A \urcorner) \urcorner $$, for any sentence form $A$ in the language of $T$, and $f(n)=0$ otherwise. And of course $\mathcal G_f$ is the arithmetical formula representing $f$ in $T$.\bigskip

Since $h$ is computable, then there would be a formula $\mathcal G_h(x,y)$ representing $h$ in theory $T$. Namely: $$\vdash_T \forall y \, (\mathcal G_h( ^o\ulcorner s \urcorner, y) \to y= \, ^o \ulcorner D(E_s(v)) \urcorner)$$ Then we define: $$ \operatorname {Dgn}(x,y) \iff \mathcal G_h(y,x)$$

When minimality is not assumed, we do the same but on the respective $E_s(v)$.

\section{References}
\begin{enumerate}
    \item Kurt Gödel, 1931, "Über formal unentscheidbare Sätze der Principia Mathematica und verwandter Systeme, I", Monatshefte für Mathematik und Physik, v. 38 n. 1, pp. 173–198. doi:10.1007/BF01700692

    \item —, 1931, "Über formal unentscheidbare Sätze der Principia Mathematica und verwandter Systeme, I", in Solomon Feferman, ed., 1986. Kurt Gödel Collected works, Vol. I. Oxford University Press, pp. 144–195. ISBN 978-0195147209. The original German with a facing English translation, preceded by an introductory note by Stephen Cole Kleene.
    \item Barkley Rosser (September 1936). "Extensions of some theorems of Gödel and Church". Journal of Symbolic Logic. 1 (3): 87–91. doi:10.2307/2269028. 
    \item \href{https://en.wikipedia.org/wiki/Rosser%27s_trick#Background}{Wikipedia. Rosser's trick. Background. 13/08/2023}
    \item \href{https://en.wikipedia.org/wiki/Diagonal_lemma#Proof}{Wikipedia. Diagonal Lemma. Proof. 13/08/2023} 
\item Pudlák, Pavel, On the length of proofs of finitistic consistency statements in first order theories, Logic colloq. ’84, Proc. Colloq., Manchester/U.K. 1984, Stud. Logic Found. Math. 120, 165-196 (1986). \href{https://zbmath.org/?q=an:0619.03037}{ZBL0619.03037}.

\href{https://citeseerx.ist.psu.edu/document?repid=rep1&type=pdf&doi=54952ff2ef51c3bea49b4dfaffd656a372f048db}{Citeseer pdf}
  \item Joel David Hamkins (https://mathoverflow.net/users/1946/joel-david-hamkins), Can this provide an example of incompleteness under the assumption of mere consistency?, URL (version: 2023-07-21): https://mathoverflow.net/q/451243
    
\end{enumerate}

\end{document}